[1994/12/01]
\documentclass[10pt, reqno]{amsart}
\usepackage[english, activeacute]{babel}
\usepackage{amsmath,amsthm,amsxtra}
\usepackage{epsfig}
\usepackage{amssymb}
\usepackage{dcolumn}
\usepackage{bm}
\usepackage{latexsym}
\usepackage{amsfonts}
\usepackage{hyperref}
\usepackage{pdftricks}
\begin{document}


\title{On the nonlinear stability of mKdV breathers}
\author{Miguel A. Alejo}
\address{Department of Mathematical Sciences, University of Copenhagen,  Denmark}
\email{Miguel.Alejo@math.ku.dk}
\author{Claudio Mu\~noz}
\address{Department of Mathematics, The University of Chicago, Chicago IL, USA}
\email{cmunoz@math.uchicago.edu}
\date{September, 2012}
\thanks{}


\chardef\bslash=`\\ 
\newcommand{\ntt}{\normalfont\ttfamily}

\newcommand{\cn}[1]{{\protect\ntt\bslash#1}}
\newcommand{\pkg}[1]{{\protect\ntt#1}}
\newcommand{\fn}[1]{{\protect\ntt#1}}
\newcommand{\env}[1]{{\protect\ntt#1}}
\hfuzz1pc 


\newtheorem{thm}{Theorem}[section]
\newtheorem{cor}[thm]{Corollary}
\newtheorem{lem}[thm]{Lemma}
\newtheorem{prop}[thm]{Proposition}
\newtheorem{ax}{Axiom}
\newtheorem*{thma}{{\bf Main Theorem}}
\newtheorem*{thmap}{{\bf Theorem A'}}
\newtheorem*{thmb}{{\bf Theorem B}}
\newtheorem*{thmc}{{\bf Theorem C}}
\newtheorem*{cord}{{\bf Corolary D}}
\newtheorem{defn}[thm]{Definition}

\theoremstyle{remark}
\newtheorem{rem}{Remark}[section]
\newtheorem*{notation}{Notation}
\newtheorem{Cl}{Claim}

\numberwithin{equation}{section}

\newcommand{\secref}[1]{\S\ref{#1}}
\newcommand{\lemref}[1]{Lemma~\ref{#1}}



%

\newcommand{\R}{\mathbb{R}}
\newcommand{\N}{\mathbb{N}}
\newcommand{\Z}{\mathbb{Z}}
\newcommand{\T}{\mathbb{T}}
\newcommand{\Q}{\mathbb{Q}}
\newcommand{\Com}{\mathbb{C}}
\newcommand{\la}{\lambda}
\newcommand{\pd}{\partial}
\newcommand{\wqs}{{Q}_c}
\newcommand{\ys}{y_c}
\newcommand{\al}{\alpha}
\newcommand{\bt}{\beta}
\newcommand{\ga}{\gamma}
\newcommand{\de}{\delta}
\newcommand{\te}{\theta}
\newcommand{\arctanh}{\operatorname{arctanh}}
\newcommand{\spawn}{\operatorname{span}}
\newcommand{\sech}{\operatorname{sech}}
\newcommand{\dist}{\operatorname{dist}}
\newcommand{\re}{\operatorname{Re}}
\newcommand{\ima}{\operatorname{Im}}
\newcommand{\diam}{\operatorname{diam}}
\newcommand{\supp}{\operatorname{supp}}
\newcommand{\pv}{\operatorname{pv}}
\newcommand{\sgn}{\operatorname{sgn}}
\newcommand{\vv}[1]{\partial_x^{-1}\partial_y{#1}}

\newcommand{\Lp}{\mathcal{L}_{+,\ve}}
\newcommand{\Lm}{\mathcal{L}_{-,\ve}}
\def\bm{\left( \begin{array}{cc}}
\def\endm{\end{array}\right)}
\def\YY{{\mathcal Y}}

 \providecommand{\abs}[1]{\lvert#1 \rvert}
 \providecommand{\norm}[1]{\lVert#1 \rVert}
 \newcommand{\sublim}{\operatornamewithlimits{\longrightarrow}}
 \newcommand{\ex}{{\bf Example}:\ }
\newcommand{\ve}{\varepsilon}
\newcommand{\fin}{\hfill$\blacksquare$\vspace{1ex}}

\newcommand{\be}{\begin{equation}}
\newcommand{\ee}{\end{equation}}
\newcommand{\ba}{\begin{equation*}}
\newcommand{\ea}{\begin{equation*}}
\newcommand{\bea}{\begin{eqnarray}}
\newcommand{\eea}{\end{eqnarray}}
\newcommand{\bee}{\begin{eqnarray*}}
\newcommand{\eee}{\end{eqnarray*}}
\newcommand{\ben}{\begin{enumerate}}
\newcommand{\een}{\end{enumerate}}
\newcommand{\nonu}{\nonumber}


\newcommand{\interval}[1]{\mathinner{#1}}

\newcommand{\eval}[2][\right]{\relax
  \ifx#1\right\relax \left.\fi#2#1\rvert}

\newcommand{\envert}[1]{\left\lvert#1\right\rvert}
\let\abs=\envert

\newcommand{\enVert}[1]{\left\lVert#1\right\rVert}
\let\norm=\enVert

\begin{abstract}
\textcolor{black}{Breather modes of the mKdV equation on the real line are known to be elastic under collisions with other breathers and solitons. 
This fact indicates very strong stability properties of breathers. In this note we describe a rigorous,
 mathematical proof of the stability of breathers under a class of small perturbations. Our proof involves} the existence of a 
nonlinear equation satisfied by all breather profiles, and a new Lyapunov functional 
which controls the dynamics of small perturbations and instability modes. In order to construct such a functional, 
we work in a subspace of the energy one. However, our proof introduces new ideas in order to attack the corresponding 
stability problem in the energy space. Some remarks about the sine-Gordon case are also considered.
 \end{abstract}
\maketitle \markboth{Stability of breathers} {Miguel A. Alejo and Claudio Mu\~noz}
\renewcommand{\sectionmark}[1]{}

\section{Introduction}

\medskip

Breathers have become a paradigm of beauty and complexity in nonlinear integrable systems \cite{AC}. Although present in only 
some very particular models (sine-Gordon (SG), modified Korteweg-de Vries (mKdV), nonlinear Schr\"odinger (NLS), among others),
 their  fascinating mixed behavior, combining oscillatory and soliton character, has attracted the attention of many 
researchers in the last thirty years \cite{Sch,BMW}. From the physical point of view,  breather solutions seem to be 
relevant to localization-type phenomena in optics, condensed matter physics and biophysics \cite{Au}.  They also play 
an important role in the modeling of freak and rogue waves events on surface gravity waves and also of internal waves 
in the stratified ocean, in Josephson junctions and even in nonlinear optics. See \cite{DyTr} for a representative set of these examples.
\smallskip

\smallskip
Many authors have considered the study of dynamical properties of integrable systems. In addition to the existence 
of well-known soliton and multi-soliton solutions, and the elastic character of the corresponding interaction 
\cite{AC}, one has to add that, from the Inverse Scattering Theory, the evolution of a rapidly decaying initial
 data can be described by purely algebraic methods. \textcolor{black}{Bounded energy solutions use to decompose into a superposition
of radiation plus a  very particular set  of nonlinear elements (this is called the {\it soliton resolution conjecture},} see e.g. Schuur \cite{Sch}).
 In such a decomposition, and for a very particular class of integrable 
models, not only solitons are allowed to appear, but also breathers and even more complex solutions. The emergence 
of solitons is not a surprise, since they are stable in the energy space \cite{Benj}, and even in less regular spaces, 
where only the mass is conserved  \cite{MV}.

\smallskip

However, in order to fully understand the aforementioned soliton resolution conjecture, \textcolor{black}{ 
the stability and instability properties of breathers should be taken into account. Numerical computations 
suggest a positive answer to the question of the mKdV breather stability on the real line \cite{AGV}.} However, 
the simple problem of a rigorous proof for their orbital stability has become a \textcolor{black}{ very challenging} open problem. 

\medskip

In this note, we announce that breathers solutions of the completely integrable mKdV equation \cite{AC} 
\be\label{mKdV}
u_{t}  +  (u_{xx} + u^3)_x =0,
\ee
are globally stable under suitable perturbations \cite{AM}. In more mathematical terms, we prove that breathers
 are \emph{orbitally stable} for initial data  which is a small perturbation of any breather solution, in a more
 regular subspace of the corresponding energy space.  Our proof also works in the case of the sine-Gordon equation \cite{AC}
\be\label{SG}
u_{tt} -u_{xx}  + \sin u =0,
\ee
provided a suitable linearized problem has the adequate spectral properties.

\medskip

\section{Breathers}

\medskip

\textcolor{black}{In what follows, every integral below is taken over $\R$.}

\smallskip

First of all, and for the sake of completeness, we introduce the standard Sobolev spaces $L^2$ and $H^k$:
$L^2:= \{ u : \int u^2 <+\infty\},$ $H^1 := \{ u\in L^2: \int u_x^2 <+\infty\}$ (the energy space), and 
$H^2 := \{ u \in H^1 : \int u_{xx}^2 <+\infty\}$ and so on. On the other hand, given a solution $u(t)$ of
 (\ref{mKdV}), the conservation laws at the $H^1$-level of regularity are the \emph{mass}
\be\label{M1}
M[u](t)  :=  \frac 12 \int u^2 = M[u](0),
\ee 
and \emph{energy} 
\be\label{E1}
E[u](t)  :=  \frac 12 \int u_x^2 -\frac 14 \int u^4 = E[u](0).
\ee 
A satisfactory Cauchy theory is also present at that level of regularity, see e.g. Kenig-Ponce-Vega \cite{KPV}.
 Solutions $u(t)$ constructed with initial data $u(0) \in H^1$ are globally well-defined: one has $u(t)\in H^1$
 for all $t$, and mass and energy are conserved along the corresponding flow. 

\smallskip

Roughly speaking, breathers are time periodic, spatially localized solutions of (\ref{mKdV}).  Indeed, given 
$\al, \bt>0$,  $x_1,x_2\in \R$, the breathers of mKdV \eqref{mKdV} are  given by 
\be\label{B}
 B_{\al, \bt}(t,x;x_1,x_2)   :=  2\sqrt{2} \partial_x \Big[ \arctan \Big( \frac{\bt}{\al}\frac{\sin(\al (x+\delta t + x_1))}{\cosh(\bt (x+\ga t +x_2))}\Big) \Big],
\ee
with $\delta := \al^2 -3\bt^2$ and $\ga := 3\al^2 -\bt^2$ (see \cite{W1}). Here $\al$ and $\bt$ represent the
 two scaling parameters of the breather (note that $\al$ is also an oscillatory parameter), and $-\ga$ describes
 the velocity of the solution. Note that $B_{\al, \bt}$ is periodic in time, but not in space, and this will be
 essential in our proof. If we take the limit $\bt/\al\ll1 $ in (\ref{B}), this allows to simplify the expression 
for the breather to
\[
 B_{\al, \bt}(t,x;0,0)   \approx 2\sqrt{2}\bt \cos(\al (x+\delta t))\sech (\bt (x+\ga t)) + O\Big(\frac{\bt}{\al}\Big),
\]
and from a qualitative point of view, it shows explicitly its wave packet like nature, as an oscillation modulated
 by an exponentially decaying function. 
 
 \smallskip
 
\textcolor{black}{We are ready to state the main result in \cite{AM} (here $B_{\al,\bt}$ refers to the breather solution \eqref{B}).}
 
\textcolor{black}{\begin{thm}[$H^2$-stability of mKdV breathers]\label{T1} Let $\al, \bt >0$. There exist parameters $\eta_0, A_0$, such
 that the following holds.  Consider $u_0 \in H^2(\R)$, and assume that there exists $\eta \in (0,\eta_0)$ such that 
\be\label{In}
\|  u_0 - B_{\al,\bt}(0,\cdot;0,0) \|_{H^2(\R)} \leq \eta.
\ee
Then there exist $x_1(t), x_2(t)\in \R$ such that the solution $u(t)$ of the Cauchy problem for the mKdV equation (\ref{mKdV}), with initial data $u_0$, satisfies
\be\label{Fn1}
\sup_{t\in \R}\big\| u(t) - B_{\al,\bt}(t,\cdot; x_1(t),x_2(t)) \big\|_{H^2(\R)}\leq A_0 \eta.
\ee
\end{thm}
Finally, note that the constant $A_0$ may depend on $\al$ and $\beta$, but it is independent of $\eta.$  There are explicit bounds on the variation of the parameters $x_1(t)$ and $x_2(t)$;  see \cite{AM} for more details.
}
 \medskip

\section{Ideas of the proof}

\medskip

In order to explain the main ideas involved in the proof of Theorem \ref{T1}, let us start by recalling the proof for the case of mKdV solitons. 
First of all, solitons are regarded as minimizers of a constrained functional in the $H^1$-topology. They are given by the expression
\be\label{Sol}
u(t,x) = Q_c (x-ct), \quad Q_c(s) := \sqrt{c} Q(\sqrt{c} s), \quad c>0,
\ee
with $Q (s):= \frac{\sqrt{2}}{\cosh (s)} =2\sqrt{2} \partial_s[\arctan(e^{s})].$ By replacing (\ref{Sol})
 in (\ref{mKdV}), one has that $Q_c>0$ satisfies the nonlinear ODE
\be\label{ecQc}
Q_c'' -c\, Q_c +Q_c^3=0, \quad Q_c\in H^1.
\ee
A very useful method to guess either the stable or unstable character of a soliton is given by the following positivity criteria:
 since the mass of a soliton is given by $M[Q_c]= 2\sqrt{c}$, then
\[
\partial_c M[Q_c] = c^{-1/2}>0.
\]
Whenever the expression above is either zero or negative, the soliton is actually unstable \cite{Benj,MMgafa}. With 
this information in mind, one of the main ingredients of the arguments employed in some of the results mentioned 
above is the introduction of a suitable \emph{Lyapunov functional}, \emph{invariant in time} and such that the 
soliton is a corresponding \emph{extremal point}. 
For the mKdV case, this functional is given by
\[
H[u](t) = E[u](t) + c \, M[u](t),
\]
where $c>0$ is the scaling of the solitary wave, and $E[u]$, $M[u]$ are given in (\ref{M1})-(\ref{E1}). A simple
 computation shows that for any $z(t)\in H^1$ small,
\bea\label{Expa1}
H[Q_c+z](t) & =&  H[Q_c] \textcolor{black}{-} \int z(Q_c''-cQ_c +Q_c^3)  + \mathcal Q\textcolor{black}{[z]}(t)  + O(\|z(t)\|_{H^1}^3).
\eea
The zero order term above is independent of time, while the first order term is zero from (\ref{ecQc}). \textcolor{black}{It turns 
out that the second order term \[
\mathcal Q[z](t)=\frac{1}{2}\int [z_x^2 + (c-3Q^2)z^2]
\]}
 is positive definite modulo two directions\textcolor{black}{\footnote{A quadratic form $Q$ 
has a negative direction $z \neq 0$ if $Q[z]<0$. Note that $\la z$ is also a negative direction, for all real-valued  $\la\neq 0$.}}, related to the invariance 
of the equation under shift and scaling transformations. \textcolor{black}{Time-dependent} parameters are then introduced in order to remove
 these unstable modes. Once these directions are controlled, the stability property follows from (\ref{Expa1}).

\smallskip

In \cite{AM}, we follow similar lines to prove the $H^2$-stability theorem of mKdV breathers \ref{T1}. First of all, recall the mass and energy of breathers:
\be\label{ME}
M[B_{\al,\bt}] =4\bt, \quad E[B_{\al,\bt}] =  \frac 43\bt \ga.
\ee
From the first identity, standard stability tests can be deduced:
\[
\partial_\al M[B_{\al,\bt}] =0, \quad \partial_\bt M[B_{\al,\bt}] =4>0,
\]
which suggest a kind of critical character of the scaling parameter $\al$, and the subcritical behavior of $\bt$,
 as expected. However, we will show in the next lines that, with enough regularity on hand, $\al$ behaves as a 
stable direction.

\smallskip

The first step into that direction is the following. \textcolor{black}{We have found an explicit, \emph{nonlinear elliptic equation} 
satisfied by any breather $B=B_{\al,\bt}$, for any fixed $t\in \R$, which is the following:  }
\be\label{EcB}
B_{(4x)} -2(\bt^2 -\al^2) (B_{xx} + B^3)  +(\al^2 +\bt^2)^2 B+ 5 BB_x^2 + 5B^2 B_{xx} + \frac 32 B^5 =0.
\ee
The proof of this result is involved and requires several new identities. In particular, a key step for the proof 
is to establish the simpler second order identity
\[
 B_{xt} +  2B  \int_{-\infty}^x\!\!\! BB_t  = (\al^2 +\bt^2)^2 B+2(\bt^2 -\al^2) \int_{-\infty}^x\!\!\! B_t.
\]
It appears that (\ref{EcB}) cannot be directly obtained from the original ideas of Lax \cite{LAX1}, since breathers 
do not decouple into solitons as time evolves.

\smallskip

\textcolor{black}{With the identity (\ref{EcB}) on hand, one can compute explicit expressions for the linear operator associated to it.}
 Indeed, if $B_1 := \partial_{x_1} B_{\al,\bt}$ and $B_2 := \partial_{x_2} B_{\al,\bt}$, 
$$
\mathcal L[B_1] =\mathcal L[B_2] =0,
$$
where $\mathcal L$ is the unbounded, self-adjoint operator
\bea\label{L1}
\mathcal L [z](x;t) &  := & z_{(4x)}(x) -2(\bt^2 -\al^2) z_{xx}(x) +(\al^2 +\bt^2)^2 z(x)  + 5B^2 z_{xx}(x)  \nonu \\
& &  + 10BB_x z_x(x)  + \ \big[ 5B_x^2  +10 BB_{xx}  + \frac {15}2B^4 -6(\bt^2 -\al^2) B^2 \big] z(x),\nonu\\
& &
\eea
defined in the dense subspace $H^4$ of $L^2$. Additionally, if $\Lambda_\al B:= \partial_\al B$ and $\Lambda_\bt B :=\partial_\bt B$, one has
\bea
\mathcal L[\Lambda_\al B] & =&  -4\al [ B_{xx} + B^3 +(\al^2+\bt^2)B], \label{LLA} \\
 \mathcal L[\Lambda_\al B] & =&  4\bt [ B_{xx} + B^3 - (\al^2+\bt^2)B].  \label{LLB}
\eea
Finally, if we define $B_0 :=\frac{\al\Lambda_{\bt} B + \bt \Lambda_\al B}{8\al\bt (\al^2 +\bt^2)}$, we get 
\be\label{B0}
\mathcal L[B_0] =-B.
\ee
These expressions will be useful in last part of the proof.

\smallskip

A detailed study of the linear operator $\mathcal L$ (\cite{AM}) reveals that, for any $t\in\R$ its kernel is 
isolated and it is spanned by $B_1$ and $B_2$, it has continuous spectrum given by the interval $[(\al^2+\bt^2)^2,+\infty)$ 
in the case $\beta\geq \al$, and $[ 4\al^2 \bt^2 ,+\infty)$ in the case $\beta< \al$. Even more surprising, it has just
 one negative eigenvalue. The first hint of this result can be obtained by computing the quadratic form associated to
 $\mathcal L$ for the elements which are the first candidates to be negative directions, namely the scaling parameters 
$\Lambda_\al B$ and $\Lambda_\bt B $. Indeed, from (\ref{LLA}) and (\ref{ME}) one has
\[
\int \mathcal L[\Lambda_\al B]\Lambda_\al B  = 4\al \partial_\al E[B] = 32\al^2 \bt>0,
\]
\begin{align*}
\int \mathcal L[\Lambda_\bt B]\Lambda_\bt B & =  -4\bt \partial_\bt E[B] - 16\bt (\al^2 +\bt^2) =  -32 \al^2 \bt <0.
\end{align*}
\textcolor{black}{A rigorous proof of the existence of a unique negative eigenvalue makes use of the spectral theory developed by  L. Greenberg \cite{Gr}.} 
As a conclusion, it is expected that, as in the standard soliton case, the scaling parameter $\beta$ describes an 
instability mode, which need to be removed, or well-estimated.

\smallskip

Finally, the more important ingredient appearing in the proof of \ref{T1} in \cite{AM} is the introduction of a Lyapunov functional, well-defined in
 the $H^2$ topology, and for which breathers are surprisingly not only extremal points, but also local minimizers, 
up to symmetries. This functional is a suitable combination of the mass and energy \eqref{M1}-\eqref{E1}, and a 
third conserved  quantity, defined at the $H^2$-level by
\be\label{F1}
F[u](t)    :=    \frac 12 \int u_{xx}^2 -\frac 52 \int u^2u_x^2  + \frac 14 \int u^6.
\ee
(The Cauchy theory for $H^2$ initial data is well understood and $F$ remains constant along the flow.) We define 
this new Lyapunov functional as follows. Let $B=B_{\al,\bt}$ be any breather with scalings $\al$ and $\bt$, and
 $M[u]$ and $E[u]$ introduced in (\ref{M1}), (\ref{E1}). We denote, for any time $t$,
\be\label{H1}
\mathcal{H}[u] :=  F[u] + 2(\bt^2-\al^2) E[u] + (\al^2 +\bt^2)^2 M[u].
\ee
This functional is reminiscent of the one appearing in the foundational paper by Lax \cite{LAX1}, concerning the 
 2-soliton solution of the Korteweg-de Vries (KdV) equation, and generalized to the KdV $N$-soliton states by 
Maddocks-Sachs \cite{Benj}. This idea has been successfully applied to several 2-soliton problems, for which the 
dynamics decouples into well-separated solitons as time evolves, see e.g. \cite{HPZ}. However, there was no 
evidence that this technique could be generalized to the case of even more complex solutions, such as breathers. 
Compared with those results, our proofs are involved, since there is no mass splitting as $t\to +\infty$. 

\smallskip

Coming back to the proof, it is clear that $\mathcal H[u]$ is a conserved quantity, well-defined for $H^2$-solutions 
of (\ref{mKdV}). Moreover, one has the following decomposition: let $z\in H^2$ be any function with sufficiently small
 norm, and $B=B_{\al,\bt}$ be any breather solution.  Then, for all $t\in \R$,  one has
\be\label{EE}
\mathcal{H}[B+z] - \mathcal{H}[B]  = \frac 12\mathcal Q[z] + \mathcal N[z],
\ee
with $\mathcal Q$ being the quadratic form defined by $\mathcal L$, namely $\mathcal Q[z] := \int z\mathcal L[z]$,
 and $\mathcal N[z]$ satisfies $|\mathcal N[z] | \lesssim \|z\|_{H^2}^3.$ Note that in this result, \textcolor{black}{the linear
 term in $z$ \emph{vanishes} since it is proportional to the left hand side of (\ref{EcB}).} 

\smallskip

From the previous analysis, we see that the functional \eqref{H1} will allow us to control the dynamics of the
 perturbative term $z(t)$, provided we manage the three instability directions that appear as a consequence of
the symmetries satisfied by \eqref{mKdV}: the two dimensional kernel $B_1,B_2$, and the eigenfunction corresponding
 to the negative  eigenvalue of $\mathcal L$. We \textcolor{black}{ vary} in time the parameters $x_1$ and $x_2$ in (\ref{B}) to 
satisfy, for all time, the orthogonality conditions
\be\label{K}
\int B_1 z = \int B_2 z =0, 
\ee
with $z(t) := u(t) -B_{\al,\bt}(t; x_1(t),x_2(t))$. This is an absolutely necessary condition in order to obtain 
an orbital stability property. However, we do not modulate the scaling instability $\bt$. Instead, we control the
 dynamics by first replacing the corresponding negative mode by a more tractable direction, the breather itself 
(this technique was first introduced by Weinstein in \cite{We1}): one has
\be\label{Coer}
\| z(0)\|_{H^2}^2 + \|z(t)\|_{H^2}^3 \gtrsim \mathcal Q[z(t)] \gtrsim \|z(t)\|_{H^2}^2 -\abs{\int \!\! zB(t)}^2,
\ee
provided (\ref{K}) is satisfied. Roughly speaking, this property is consequence of (\ref{B0}) and the fact that
 $\int B_0B\neq 0$, that is, $B_0$ and $B$ are not orthogonal, and $B_0$ is also a negative direction. Finally, 
we use the mass conservation law to give a new estimate on the last term in (\ref{Coer}), which is at the first 
order independent of time, and only depending on $z(0)$.
This last element of the proof  allows to control the term $z(t)$ for all times, proving the stability property. 

\medskip

\section{Discussion}

\medskip

Our results emphasize, and demonstrate, at the rigorous level, some deep connections between breathers and 
the 2-soliton solution of mKdV. This connection comes from the structure of the Lyapunov functional \eqref{H1} that 
we have used as an essential tool in the proof. The main theorem is stated for an $H^2$ perturbation, although we 
believe  that it can be improved to reach the $H^1$ level, but with a harder proof. \textcolor{black}{ We also believe that our proof
 gives some insights to approach the study of the instability of \emph{periodic} arrays of breathers, 
as numerically stated in \cite{KKSH,KKS}. In order to prove this fact, it is necessary to find a suitable Lyapunov functional for which lattice breathers are critical points. }

\smallskip

Concerning the SG case, our proof also applies with no significant modifications (treating SG as a matrix operator
 problem for $(u,u_t) \in H^2 \times H^1$), provided a suitable linearized operator (of matrix type and of fourth 
order in terms of derivatives), has only one negative eigenvalue, and the associated kernel is nondegenerate. For 
the moment, a proof of that result has escaped to us. Therefore, if these conditions are satisfied, SG breathers are 
$H^2\times H^1$ stable.  

\smallskip

\textcolor{black}{However, it is highly expected that breathers could not survive under nontrivial perturbations of the equation, as is showed in \cite{BMW}, in the SG case. }

\medskip

\section{Conclusion}

\medskip

 In this note we have presented some details of the proof of the nonlinear stability of breathers for the mKdV equation.
 Our conclusion is that mKdV breathers are nonlinear stable at the $H^2$ level of regularity. No systematic work on this
 kind was reported before, probably as a direct consequence of the difficulties that arise by the complexity of the
 breather solution and the computations involved. As we explained above, the introduction of the Lyapunov functional
 \eqref{H1} for which breathers are local minimizers, up to symmetries, is the essential point. This functional 
 allows to control the perturbative terms and the instability directions that appear during of the dynamics. 
Finally, we indicate that  our arguments are general and, in principle, can be applied to several equations with 
breather solutions.  We expect to consider some of these problems in a forthcoming publication.

\smallskip

{\bf Acknowledgments}. \textcolor{black}{We wish to express our sincere thanks to the anonymous referees for their careful reading and useful 
 suggestions.} We would also like to thank Y. Martel, F. Merle and L. Vega for useful comments and discussions on this work.

\end{document}